\documentclass{amsart}
\usepackage{amsmath, amssymb, amsthm}
\usepackage[utf8]{inputenc}
\usepackage[T1]{fontenc}
\usepackage[colorlinks=true, allcolors=blue]{hyperref}
\usepackage{float}

\theoremstyle{plain}
\newtheorem{theorem}{Theorem}[section]

\newtheorem{corollary}[theorem]{Corollary}
\newtheorem{proposition}[theorem]{Proposition}

\title[Jordan nilpotent group rings]{On Jordan nilpotency of index n in group rings}
\author[Cristo]{Osnel Broche Cristo}\address{Departament of Mathematics and Applied Mathematics  \\ Federal University of Lavras \\   37200-000, Lavras - MG, Brazil}
\email{osnel@dex.ufla.br}
\author[Hoyos]{Mariana Garabini Cornelissen Hoyos}\address{Departament of Estatistic, Physics and  Mathematics \\ Federal University of São João Del Rei \\ Ouro Branco - MG, Brazil}
\email{mariana@ufsj.edu.br}
\author[Rodrigues]{Rodrigo Lucas Rodrigues}\address{Department of Mathematics\\ Federal University of Cear\'a\\ Pici Campus, Block 914\\ 60455-760, Fortaleza - CE, Brazil}
\email{rodrigo@mat.ufc.br}

\begin{document}

\begin{abstract} 
Let $RG$ be the group ring of an arbitrary group $G$ over an associative and commutative ring $R$ with identity. In this paper, we provide necessary and sufficient conditions for $RG$ to be Jordan nilpotent of index $n \geq 2$. In particular, we establish a connection between the study of Jordan nilpotency in group rings and the study of Lie nilpotency in these algebras. We also classify the Jordan nilpotency of $RG$, over rings of arbitrary characteristic, up to the index $6$.
\end{abstract}
	
	\maketitle



\section{Introduction}

An algebra is said to be Lie if it satisfies the identities $x^2 = 0$ (anticommutativity) and $x(yz) + y(zx) + z(xy) = 0$ (Jacobi identity). An algebra is called Jordan if it satisfies the identities $xy = yx$ (commutativity) and $(x^2 y) x = x^2 (yx)$ (Jordan identity). 

Given an associative algebra $A$, considering the same underlying vector space as $A$, we can construct the algebra $(A, [\cdot, \cdot])$ with the commutator $[x, y] = xy - yx$ as the multiplication and the algebra $(A, \circ)$ with the Jordan product $x \circ y = xy + yx$. It is well known that $(A, [\cdot, \cdot])$ is a Lie algebra and $(A, \circ)$ is a Jordan algebra.  

We say that $(A, [\cdot,\cdot])$ is a Lie nilpotent algebra if, for some positive integer $n \geq 2$, $ [[ \cdots [[x_1, x_2],x_3] \cdots ],x_n] = 0$ holds for all $x_1, x_2, \ldots, x_n \in A$. Similarly, $(A, \circ)$ is a Jordan nilpotent algebra if, for some positive integer $n \geq 2$, $(\cdots ((x_1 \circ x_2) \circ x_3) \cdots ) \circ x_n =0$ holds for all $x_1, x_2, x_3, \ldots, x_n \in A$. The smallest $n$ for which a Lie (Jordan) algebra is nilpotent is called the index of nilpotency.  

One of the interesting problems which arose was studying the Lie properties of group algebras, particularly determining necessary and sufficient conditions for a group algebra to be Lie nilpotent. In 1973, Passi et al. \cite{pps} proved that for a field $K$,  $KG$ is Lie nilpotent if and only if either the characteristic of $K$ is zero and $G$ is an abelian group or characteristic of $K$ is a prime $p$, $G$ is nilpotent and $G'$ is a finite $p$-group. In 1986, Levin and Rosenberger \cite{lr} proved that $RG$ is Lie nilpotent of index up to $4$ and ${\rm char}(R) = 2$ if and only if $G$ is a nilpotent group of class at most 2 and $G'$ has exponent 2 and order 2 or 4.

Let $K$ be a field of characteristic $p > 0$. In 2014, Chandra and Sahai \cite{cs} classified group algebras $KG$ which are Lie nilpotent of index $k \leq 6$.  In 2019, Sahai and Sharan \cite{ss} gave a description of Lie nilpotent group algebras $KG$ of index $k = 7$ or $k = 8$. Since then, there is no known classification for Lie nilpotency of group algebras for higher indices with such an approach.

In a natural way it is also of interest to investigate Jordan nilpotent group algebras. In this direction, Goodaire and Polcino Milies \cite{gpm} characterized the group rings $RG$ satisfying $(x_1 \circ x_2) \circ x_3 = 0$ for all $x_1, x_2, x_3 \in RG$, and consequently, Jordan nilpotent group algebras $RG$ of index $\leq 3$ are classified. The classification of group rings that are Jordan nilpotent of index $4$ was established by the authors in \cite{n4}.

In this article, we wish to investigate the Jordan nilpotency of index $n$ for $RG$. Throughout this paper, $R$ is an associative and  commutative ring with identity, $G$ is an arbitrary group and $RG$ is the group ring of $G$ over $R$. The characteristic of $R$ will be denoted by ${\rm char}(R)$, $C_k^{j}$ will denote the direct product of $j$ copies of $C_k$, where $C_k$ is the cyclic group of order $k$, and $G' = \left\langle (g, h) \ | \ g, h \in G \right\rangle$ stands for the derived subgroup of $G$, where $(g,h) = g^{-1}h^{-1}gh$. We can inductively define the derived series of subgroups by setting $\gamma_{1}(G) = G$ and $\gamma_{k}(G) = (\gamma_{k-1}(G), G)$, for $k \geq 2$. In particular, $\gamma_3(G) = 1$ if and only if $\gamma_2(G) \subseteq \mathcal{Z}(G)$, where $\mathcal{Z}(G)$ is the center of $G$.  

The remainder of the paper is organized as follows. Section 2 is devoted to some basic and essential results for the understanding and development of this work. Section 3 deals with the case of characteristic 2, while Section 4 provides necessary and sufficient conditions for $RG$ to be Jordan nilpotent of index $n \geq 2$ and Section 5 classifies Jordan nilpotent group rings of index $n \leq 6$ and $char(R) \not= 2$.

\section{Some Basic Facts}

We know that the Jordan nilpotency of the group ring $RG$ depends on both $G$ and $R$, but primarily on the structure of $G$. We begin this section by showing that, with respect to $R$, it depends only on the characteristic of $R$.


Let $R_1$ and $R_2$ be rings such that
${\rm char}(R_1)={\rm char}(R_2)=k$.
Then
\[
\underbrace{(\cdots (R_1G\circ R_1G)\cdots \circ R_1G)\circ R_1G}_{n\ \mathrm{times}}=0
\Leftrightarrow 
\underbrace{(\cdots (G\circ G)\cdots \circ G)\circ G}_{n\ \mathrm{times}}=0.
\]
Equivalently,
\[
\underbrace{(\cdots (G\circ G)\cdots \circ G)\circ G}_{n\ \mathrm{times}}
\subseteq kR_1G.
\]
Since ${\rm char}(R_1)={\rm char}(R_2)=k$, this is equivalent to
\[
\underbrace{(\cdots (G\circ G)\cdots \circ G)\circ G}_{n\ \mathrm{times}}
\subseteq kR_2G.
\]
Hence
\[
\underbrace{(\cdots (G\circ G)\cdots \circ G)\circ G}_{n\ \mathrm{times}}=0
\quad \text{in } R_2G,
\]
which is equivalent to
\[
\underbrace{(\cdots (R_2G\circ R_2G)\cdots \circ R_2G)\circ R_2G}_{n\ \mathrm{times}}=0.
\]

We have thus proved the following result.

\begin{proposition}\label{ring}
If $R_1$ and $R_2$ are rings such that ${\rm char}(R_1) = {\rm char}(R_2) = k$, then $R_1 G$ is Jordan nilpotent of index $n$ if and only if $R_2 G$ is Jordan nilpotent of index $n$.  
\end{proposition}

The result below is a well-known fact. Since we were unable to trace its origin, we decided to include its proof here for the sake of completeness.

\begin{proposition}\label{squares}
If $G' = \langle s \rangle  \cong C_2$, then $G'$ is central. Furthermore, if the square of every element of $G$ belongs to $\mathcal{Z}(G)$, then $\gamma_3(G) = 1$ and $G'$ is an elementary abelian $2$-group.
\end{proposition}
\begin{proof}
If $sg \neq gs$, for some $g \in G$, then $(g,s)=s$ and $s=1$, a contradiction. Now, suppose that the square of every element of $G$ belongs to $\mathcal{Z}(G)$. A trivial verification shows that $(x,y) = x^{-2}(xy^{-1})^{-2}y^2$, which implies $\gamma_3(G) = 1$. Furthermore, if $(g, h) = s  \not=1$, then $s^2 = 1$, since
\[h^2g = gh^2 = (gh)h = (hgs)h = h(gh)s = h(hgs)s = h^2gs^2. \qedhere \]    
\end{proof}

In what follows, we investigate when the squares of elements of $G$ lie in the center. Observe that if $RG$ is Jordan nilpotent of index $n$ then
\[ 0 = (\cdots ((\underbrace{1 \circ 1)\circ 1) \circ \cdots 1) \circ 1}_{\mbox{$n$ elements 1's}}= 2^{n-1}\] which implies that char($R$) divides $2^{n-1}$, that is, char($R$) = $2^k$, $1 \leq k \leq n-1$.  

\begin{proposition}\label{squares1}
    If $RG$ is Jordan nilpotent of index $n$, then the square of every element of $G$ belongs to $\mathcal{Z}(G)$ in the following cases:
    \begin{itemize}
    \item[(i)] $n = 2, 3$ and $4$;
    \item[(ii)] $n = 5, 6$ and ${\rm char}(R) \not=2$.
    \end{itemize}
\end{proposition}
\begin{proof}
For $n = 2,3$ and $4$, the result follows immediately from Proposition 2.2 of \cite{n4}.  

Now suppose $n=5$ and ${\rm char}(R) \neq 2$. In this case, ${\rm char}(R) \in \left\{4, 8, 16 \right\}$ and $0 = (((g \circ h) \circ g) \circ g) \circ g = g^4 h + 4g^3 h g + 6 g^2 h g^2 + 4 g h g^3 + hg^4$ clearly force $g^2 h g^2 \in \left\{g^4 h, g^3 h g, g h g^3, hg^4 \right\}$. Analyzing each possibility, we conclude that $hg^2 = g^2 h$ or $hg = gh$. Then the squares are central. 

Finally, for $n = 6$, we know that ${\rm char}(R) \ | \ 32$ and
$$0=((((g \circ h) \circ g) \circ g) \circ g) \circ g = g^5 h + 5 g^4 hg + 10 g^3 h g^2 + 10 g^2 h g^3 + 5 ghg^4 + hg^5.$$
If ${\rm char}(R) \neq 2$, then $g^3 h g^2 \in \left\{g^5h, g^4 h g, g^2 h g^3, ghg^4, hg^5\right\}$. Now, $g^3hg^2 = hg^5$ implies that $g^3 h g^2 \in \left\{g^5h, g^4 h g, g^2 h g^3, ghg^4\right\}$. In any case, we conclude that $gh^2 = h^2g$ or $gh = hg$, that is, the squares are central.  
\end{proof}

As a consequence of Propositions \ref{squares} and \ref{squares1}, we can assert the following result.

\begin{corollary}\label{commutator}
If $RG$ is Jordan nilpotent of index $n$, then $\gamma_3(G) = 1$ and $G'$ is an elementary abelian $2$-group in the following cases:
    \begin{itemize}
    \item[(i)] $n = 2, 3$ and $4$;
    \item[(ii)] $n = 5, 6$ and ${\rm char}(R) \not=2$.
    \end{itemize}   
\end{corollary}

\section{Characteristic 2}\label{Section2}
We first consider the case ${\rm char}(R) = 2$. In this case, $RG$ is Jordan nilpotent of index $n \geq 2$ if and only if $RG$ is Lie nilpotent of index $n$.



As a consequence of Proposition \ref{ring}, Theorem 5.1 in \cite{n4}, Theorems 3.3 and 3.4 in \cite{cs} and Theorems 2.5 and 2.6 in \cite{ss}, Jordan nilpotent group rings of index $n \leq 8$ can be characterized as follows.

\begin{theorem} \label{2} Let $RG$ denote the group ring of a group $G$ over a ring $R$ such that  ${\rm char}(R) = 2$. Then 
\begin{enumerate}
    \item $RG$ is Jordan nilpotent of index $2$ if and only if $G$ is abelian;
    \item $RG$ is Jordan nilpotent of index $3$ if and only if $G' \cong C_2$ (and $\gamma_3(G) = 1)$;
    \item $RG$ is Jordan nilpotent of index $4$ if and only if $G' \cong C_2^2$ and $\gamma_3(G) = 1$.
    \item $RG$ is Jordan nilpotent of index $5$ if and only if one of the following conditions holds:
\begin{itemize}
\item[(i)] $G' \cong C_4$ and $\gamma_4(G) = 1$;
\item[(ii)] $G' \cong C_2^2$, $\gamma_3(G) \neq 1$ and $\gamma_4(G) = 1$;
\item[(iii)] $G' \cong C_2^3$ and $\gamma_3(G) = 1$.
\end{itemize}
    \item $RG$ is Jordan nilpotent of index $6$ if and only if one of the following conditions holds:
\begin{itemize}
\item[(i)] $G' \cong C_2^3$, $\gamma_3(G)$ cyclic and $\gamma_4(G) = 1$;
\item[(ii)] $G' \cong C_4 \times C_2$, $G'^2 = \gamma_3(G)$ and $\gamma_4(G) = 1$;
\item[(iii)] $G' \cong C_2^4$ and $\gamma_3(G) = 1$.
\end{itemize}

\item $RG$ is Jordan nilpotent of index $7$ if and only if one of the following conditions holds:
\begin{itemize}
\item[(i)] $G' \cong C_2^5$ and $\gamma_3(G)=1$;
\item[(ii)] $G' \cong C_4 \times C_2 \times C_2$ and $\gamma_3(G) \subseteq G'^2$;
\item[(iii)] $G' \cong C_2^4$ and $\gamma_3(G) \cong C_2$;
\item[(iv)] $G' \cong C_4 \times C_2$, $G'^2 \subseteq \gamma_3(G) \cong C_2 \times C_2$ and $\gamma_4(G) = 1$;
\item[(v)] $G' \cong C_4 \times C_2$, $\gamma_3(G) \cong C_2$ and $G'^2 \cap \gamma_3(G) = 1$;
\item[(vi)] $G' \cong C_2^3$, $\gamma_3(G) \cong C_2 \times C_2$ and $\gamma_4(G) = 1$.
\end{itemize}

\item $RG$ is Jordan nilpotent of index $8$ if and only if one of the following conditions holds:
\begin{itemize}
\item[(i)] $G' \cong C_2^6$ and $\gamma_3(G) = 1$;
\item[(ii)] $G' \cong C_4 \times C_2^3$ and $\gamma_3(G) \subseteq G'^2$;
\item[(iii)] $G' \cong C_4 \times C_4$, $\gamma_3(G) \subseteq G'^2$ and $\gamma_4(G)=1$;
\item[(iv)] $G' \cong C_2^5$ and $\gamma_3(G) \cong C_2$;
\item[(v)] $G' \cong C_4 \times C_2 \times C_2, G'^2 \subseteq \gamma_3(G) \cong C_2 \times C_2$ and $\gamma_4(G)=1$;
\item[(vi)] $G' \cong C_4 \times C_2 \times C_2, \gamma_3(G) \cong C_2$ and $G'^2 \cap \gamma_3(G) = 1$;
\item[(vii)] $G' \cong C_2^4, \gamma_3(G) \cong C_2 \times C_2$ and $\gamma_4(G)=1$;
\item[(viii)] $G' \cong C_4 \times C_2, \gamma_3(G) \cong C_2 \times C_2$ and $\gamma_4(G) \cong C_2$;
\item[(ix)] $G' \cong C_2^3, \gamma_3(G) \cong C_2 \times C_2$ and $\gamma_4(G) \cong C_2$.
\end{itemize}
\end{enumerate}
\end{theorem}


\section{Necessary and sufficient conditions}

The next useful result, of reducing the Jordan nilpotency case to the Lie nilpotency case, establishes necessary conditions for $RG$ to be Jordan nilpotent of index $n \geq 2$ and will be of fundamental importance throughout the article.

\begin{theorem}\label{connection}
 If $RG$ is Jordan nilpotent of index $n \geq 2$, then for some $1 \leq i \leq n-1$, ${\rm char}(R) = 2^{n-i}$ and $LG$ is Lie  nilpotent of index at most $i+1$, for every commutative ring $L$ with identity such that ${\rm char}(L) = 2$.
\end{theorem}

\begin{proof}
Recall that if $RG$ is Jordan nilpotent of index $n \geq 2$, then ${\rm char}(R)$ divides $2^{n-1}$. Let ${\rm char}(R) = 2^{n-j}$ for some $j \in \left\{1, \ldots, n-1\right\}$. Now,

$$( \cdots ((( \cdots \underbrace{(RG \circ RG) \cdots \circ RG) \circ RG)}_{j + 1 \ {\rm times}} \circ \underbrace{1) \cdots) \circ 1}_{n - j - 1 \ {\rm times}} = 0.$$

$$2^{n - j - 1}\underbrace{(( \cdots (RG \circ RG) \cdots \circ RG) \circ RG)}_{j + 1 \ {\rm times}} = 0.$$

In particular,
$2^{n - j - 1}\underbrace{(( \cdots (G \circ G) \cdots \circ G) \circ G)}_{j + 1 \ {\rm times}} = 0$ holds. Since $G$ is a basis of $RG$ over $R$, we see that
$\underbrace{( \cdots (G \circ G) \cdots \circ G) \circ G}_{j + 1 \ {\rm times}} \subseteq 2RG$. Therefore, for every commutative ring $L$ with identity such that ${\rm char}(L) = 2$, it may conclude that
$$\underbrace{( \cdots (G \circ G) \cdots \circ G) \circ G}_{j + 1 \ {\rm times}} =0 \ {\rm in} \ LG.$$

Consequently,
 $\underbrace{( \cdots (LG \circ LG) \cdots \circ LG) \circ LG}_{j + 1 \ {\rm times}} =0,$
that is, $LG$ is Jordan nilpotent of index at most $j+1$.   
\end{proof}

The following result provides further necessary conditions for $RG$ to be Jordan nilpotent.

\begin{theorem} \label{CN}
    If $RG$ is Jordan nilpotent of index $n$ and ${\rm char}(R) = 2^{n-i}$, with $1 \leq i \leq 3$ then $G' \cong C_{2}^{i-1}$ and  $\gamma_{3}(G) = 1$.
\end{theorem}

\begin{proof} Since $RG$ is Jordan nilpotent of index $n$, we have ${\rm char}(R) = 2^{n-i}$ for some $1 \leq i \leq n-1$. Assume $i \in \{1,2,3\}$. Let $L$ be a commutative ring with identity such that ${\rm char}(L) = 2$. If $RG$ is Jordan nilpotent of index $n \geq 2$, it follows from Theorem \ref{connection}, that $LG$ is Lie nilpotent of index $j$, where $2 \leq j \leq i+1 \leq 4$. 
For $i = 1$, $G$ is abelian and the assertion follows.
Suppose now that $i \neq 1$ and $j < i+1$, then $G' \cong C_2^{j-2}$ by Theorem \ref{2}.  
Hence, $RG$ is Jordan nilpotent of index $n - (i-j+1) < n$, a contradiction. Therefore, $j = i+1$, and so $LG$ is Lie nilpotent of index $i+1$. Consequently, by Theorem \ref{2}, $G' \cong C_2^{i-1}$ and $\gamma_3(G) = 1$. 
\end{proof}

We now establish sufficient conditions for $RG$ to be Jordan nilpotent.

\begin{theorem}\label{CS}
If ${\rm char}(R) = 2^{n-j}$, $1 \leq j \leq 8$, $\gamma_3(G) = 1$ and $G' \cong C_2^{j-1}$, then $RG$ is Jordan nilpotent of index $n$.
\end{theorem}

\begin{proof} 
A standard calculation shows that for $g_i \in G$, $1 \leq k \leq j$ and $c_i \in G'$, we have:
\[ ( \cdots ( ( g_1 \circ g_2) \circ g_3) \cdots \circ g_{n-1} ) \circ g_n = 2^{n-k} \prod_{i=1}^{k-1}(1+c_i) \prod_{i=1}^{n} g_i  \] 
Since $char(R) = 2^{n-j}$, the above expression is zero, which implies that $RG$ is Jordan nilpotent of index at most $n$. 

Suppose that $RG$ is Jordan nilpotent of index $m$, with $m<n$. Then ${\rm char}(R) = 2^{m-i}$ with $1 \leq i \leq m-1$. Since ${\rm char}(R) = 2^{n-j}$, it follows that $1\leq i = m-n+j$ and moreover $i<j$ because $m<n$. By Theorem \ref{connection}, $LG$ is Jordan nilpotent of index $k$, where ${\rm char}(L)=2$ and $k \leq i+1 \leq j$. Since $k \leq 8$, Theorem \ref{2} implies that $G'$ has at most $k-2$ copies of $C_2$, contradicting the assumption that $G' \cong C_2^{j-1}$. 
\end{proof}


Gathering the results established above, we obtain the following necessary and sufficient conditions for $RG$ to be Jordan nilpotent.

\begin{theorem}\label{CNS}
Suppose ${\rm char}(R) = 2^{n-i}$ with $1 \leq i \leq 3$.  $RG$ is Jordan nilpotent of index $n$  if and only if  $G' \cong C_{2}^{i-1}$ and  $\gamma_{3}(G) = 1$.
\end{theorem}

\section{Jordan Nilpotency of index at most 6}

We are now in a position to present necessary and sufficient conditions for $RG$ to be Jordan nilpotent of index at most 6, based on the preceding results. We recall that the case where $char(R) = 2$ was addressed in Section 2. Therefore, throughout this discussion, we assume that $char(R) \not=2$.

\begin{theorem}
Let $RG$ denote the group ring of a group $G$ over a coefficient ring $R$ with ${\rm char}(R) \neq 2$. $RG$ is Jordan nilpotent of index $n$,  $2 \leq n \leq 6$, if and only if ${\rm char}(R) = 2^{n-i}$, $G' \cong C_2^{i-1}$, $1 \leq i \leq n-1$, and $\gamma_3(G) = 1$.
\end{theorem}

\begin{proof}

According to Theorem \ref{CS}, we only have to prove that the necessary conditions are true. To this end, suppose that $RG$ is Jordan nilpotent of index $n \leq 6$. Then, by Theorem \ref{connection}, ${\rm char}(R) = 2^{n-i}$, for some $1 \leq i \leq n-1$ and, by Corollary \ref{commutator}, we have that $\gamma_3(G) = 1$. It remains to show that $G' \cong C_2^{i-1}$, $1 \leq i \leq n-1$. By Theorem \ref{CN}, we just need to prove for $n=6$ and $i=4$.

In this case, we have ${\rm char}(R) = 4$ and, by Theorem \ref{connection}, $LG$ is Lie nilpotent of index at most 5, for every associative commutative ring $L$ with identity such that ${\rm char}(L) = 2$.
Furthermore, by Theorem \ref{2}, if $LG$ is Jordan nilpotent of index $j<5$, then $G' \cong C_2^{j-2}$, and consequently, by Proposition \ref{CS},
$RG$ is Jordan nilpotent of index at most $j+1<6$, a contradiction.
Therefore, $LG$ is Lie nilpotent of index 5 and we conclude, by 
Theorem \ref{2}, that one of the following conditions holds:
\begin{itemize}
\item[(I)] $G' \cong C_4$ and $\gamma_4(G) = 1$;
\item[(II)] $G' \cong C_2 \times C_2$, $\gamma_3(G) \neq 1$ and $\gamma_4(G) = 1$;
\item[(III)] $G' \cong C_2 \times C_2 \times C_2$ and $\gamma_3(G) = 1$.
\end{itemize}
Now, by Corollary \ref{commutator}, $\gamma_3(G) = 1$ and $G' \not\cong C_4$. Then (I) and (II) do not occur, that is,   $G' \cong C_2 \times C_2 \times C_2$ and $\gamma_3(G) = 1$.
\end{proof}

\section*{Acknowledgments}
This work was supported by FAPEMIG (Fundação de Amparo à Pesquisa do Estado de Minas Gerais) under grants APQ-00958-22 and RED-00133-21.

\end{document}